\documentclass[12pt]{article}

\usepackage{graphicx}
\usepackage[centertags]{amsmath}
\usepackage{amsfonts}
\usepackage{amssymb}
\newcommand{\pr}{\mbox{\sf P}}
\newcommand{\ex}{{\bf\sf E}}               %% expectation
\newcommand{\var}{\mbox{\sf Var}}

\def\bsk{\bigskip\noindent}

\newcommand{\bZ}{{\bf Z}}
\newcommand{\bx}{{\bf x}}               %% bold x (vector)

\newcommand{\al}{\alpha}                %%
\newcommand{\g}{\lambda}                %%
\newcommand{\imp}{\Rightarrow}           %% arrows
\newcommand{\startb}{\parindent0pt\bf}  %% start a boldface line
\newtheorem{thm}{Theorem}

\newtheorem{defn}{Definition}

\def\th{\theta}
\def\Th{\Theta}

\def\ind{\chi}

\begin{document}

%\begin{frontmatter}

\begin{titlepage}
	\centering
	{\Large \bfseries Heuristic Policies for Stochastic Knapsack Problem with Time-Varying Random Demand}
	\vspace{2cm}
	
	Yingdong Lu\\
	Mathematical Sciences Department \\ IBM T.J. Watson Research Center \\ Yorktown Heights, NY 10598, U.S.A. \\ Phone: 914-945-3738, Fax: 914-945-3434, \\ Email: yingdong@us.ibm.com.
	
	\end{titlepage}



\title{Stochastic Knapsack Problem with Time-Varying Random Demand: Structural Properties and Heuristics}


\author{Yingdong Lu \\Mathematical Sciences Department \\ IBM T.J. Watson Research Center \\ Yorktown Heights, NY 10598, U.S.A. \\ yingdong@us.ibm.com}

\begin{abstract}
In this paper, we consider the classic stochastic (dynamic) knapsack problem, a fundamental mathematical model in revenue management, with general time-varying random demand. Our main goal is to study the optimal policies,  which can be obtained by solving the dynamic programming formulated for the problem, both qualitatively and quantitatively. It is well-known that when the demand size is fixed and the demand distribution is stationary over time, the value function of the dynamic programming exhibits extremely useful first and second order monotonicity properties, which lead to monotonicity properties of the optimal policies. In this paper, we are able to verify that these results still hold even in the case that the price distributions are time-dependent. When we further relax the demand size distribution assumptions and allow them to be arbitrary, for example in random batches, we develop a scheme for using value function of alternative unit demand systems to provide bounds to the value function. These results confirm some of the basic understandings of the stochastic knapsack problem. As a natural consequence, we propose a new class of heuristic policies, two-dimensional switch-over policies, and discuss numerical procedures for optimizing the parameters of the policies. Extensive numerical experiments demonstrate that these switch-over policies can provide performances that are in a close neighborhood of the optimum statistically.  
\end{abstract}

%\begin{keyword}


%Revenue Management, Dynamic Programming
%\end{keyword}

%\end{frontmatter}


\section{Introduction}
\label{sec:intro}

\subsection{Motivations}
\label{sec:motivations}

Stochastic knapsack problem, also known as stochastic dynamic knapsack problem, has been studied as a fundamental mathematical model in the field of revenue management and dynamic pricing. The problem, in its rather simple format, captures the fundamental trade-off between accepting current reward and holding out for future, and potentially better opportunities. Evidently, this trade-off is not unique to revenue management, which explains the appearance of similar problem formulations, as well as results on related structural and computational properties, in many other areas such as queueing control and sequential decision making.

It is known that the stochastic knapsack problem is abundantly rich in structure, especially, when the demand at each time period can be treated as a single unit, and the statistical characterization of the price stays unchanged over time. Apart from the obvious first order monotonicities with respect to time and inventory, the value function, i.e. the maximum expected revenue can be collected, is concave in time, concave in inventory, and submodular in both time and inventory. These properties lead to important monotonicity properties in the optimal policies. However, when a multiple-unit demand model has to be adopted, one can not expect these second order monotonicity properties to hold. In fact, as we will demonstrate through a very simple numerical example, even when the demand sizes follow some really simple and common distributions, the monotonicity can break down immediately. However, heuristic policies that are motivated by those monotonicity properties in the unit demand case are proven to be asymptotically optimal, and demonstrated to be efficient through numerical experiments, see, e.g. \cite{LinLuYao2008}. 

In real world applications, it is very rare that the statistical characterization of the demand process remains unchanged over time. The probabilistic assumption on the demand usually reflects the practice of forecasting for business operations. Therefore, it is quite common that the demand forecast needs to accommodate features such as business trends and other time sensitive modifications. For example, demand for flight tickets, hotels and sports/entertainment events is often highly time-sensitive. In this paper, we are able to demonstrate that the effect of the time-varying feature on the monotonicity is not significant. More specifically, we are still able to obtain these monotonicity properties of the value function under really general assumptions on the demand pattern. In addition, we are able to connect some of these properties to multimodularity, an important concept in optimization, control and discrete analysis. With the reaffirmation of the rich structural properties on the unit demand problem, it is natural to ask if mathematically more tractable models can be used to approximate value functions in problems with general demand size distributions. Under mild assumptions on the size distributions, we are able to identify an upper bound that is related to a unit demand problem, and a lower bound that is related to a simple heuristic policy, which leads to an explicit estimation on the value function.

These qualitative results lay the foundation for our quantitative study of the optimal policy. The closed-loop switch-over policies, proposed and studied in \cite{LinLuYao2008},  are based on structural properties of the stochastic knapsack problem with stationary demand. Although the policies with optimized parameters are shown to be powerful in some cases, they are parameterized only by the remaining time in the planning horizon, hence, has a natural limitation in coping with the variation of the demand distributions. This shortcoming will be even more troublesome when facing time-varying demand. In this paper, we propose a class of heuristics of two-dimensional switch-over policies which will be parameterized by both remaining time and inventory. The switch-over points are decided through applying diffusion approximation techniques and detailed Brownian calculations. Since the one dimensional switch-over policy is just a special case in which the inventory parameter is fixed, this policy remains asymptotically optimal under the same regime of that developed in \cite{LinLuYao2008}. Meanwhile the efficiency of the heuristic, as well as its advantage over the one-dimensional switch-over policies, are demonstrated through extensive numerical experiments.

\subsection{Literature Review}
\label{sec:literature_review}

Stochastic knapsack problem is closely related to the problem of "Selling a fixed, finite inventory during a finite time period" reviewed in den Boer \cite{DenBoer2015}. Kincaid and Darling \cite{KincaidDarling1963} are believed to have pioneered the study of the dynamic pricing policies for the problem. They obtain some basic properties of the optimal acceptance and rejection strategy, and calculate the value function in some special cases. Gallego and Van Ryzin \cite{GallegovanRyzin1994} state the problem in a continuous framework, present the Hamilton-Jacobi-Bellman equation for the value function, and derive their convexity properties.  Various extensions are provided for different problem setups in  \cite{FengGallego1995},
\cite{FengGallego2000}, and \cite{FengXiao2000}. Similar models  appear in \cite{Prastacos1983},\cite{PapastavrouRajagopalanKleywegt1996}, \cite{KleywegtPapastavrou2001}, \cite{Lippman1975}, \cite{VanSlykeYoung2000}, \cite{LinLuYao2008} as well. Results that are similar to those presented in the current paper but with different approaches under a different setting can be found in  \cite{ZhaoZheng2000}. Meanwhile, there are also several studies presenting efficient heuristic policies such as the policies proposed in \cite{BitranMondschein95}, and parameterized heuristics discussed in \cite{LinLuYao2008}. 

In addition to the usual monotonicity, such as increasing, decreasing, convexity and submodularity, we observe that the value function actually can be treated as a multimodular function. Multimodularity is a concept that is  instrumental in  characterizing optimality in optimal control and discrete convex analysis. It is first defined in \cite{Hajek85},  then explored to study optimal control of queues and queueing networks, see, e.g. \cite{WeberStidham87} and \cite{GlassermanYao1994}.  Many of these results are summarized in \cite{AltmanGaujalHordijk2000}. New applications of multimodularity in supply chain management can be found in \cite{Zipkin2008}, \cite{HuhJanakiraman2010}, and \cite{LiYu2014}. The connection between multimodularity and discrete convexity is established in \cite{Murota2005}.
%Here, we will give the literature including \cite{AltmanGaujalHordijk2000}, \cite{Murota2005}, \cite{Zipkin2008}, \cite{HuhJanakiraman2010}.

%We need to give credits to where credit is due. But also need to point out that who really had the idea of multimodularity. the connection between multimodularity and discrete convexity is easy to make. 

%in dynamic programming, \cite{Prastacos1983} ,\cite{PapastavrouRajagopalanKleywegt1996}, \cite{KleywegtPapastavrou2001}, 
%queueing control \cite{Lippman1975}, and revenue management \cite{GallegovanRyzin1994},\cite{BitranMondschein95}, \cite{FengGallego1995},
%\cite{FengGallego2000},\cite{FengXiao2000},  \cite{VanSlykeYoung2000}, \cite{ZhaoZheng2000},\cite{LinLuYao2008},

\subsection{Organization of the Paper}

The rest of the paper will be organized as follows. In Sec. \ref{sec:models}, we will present the basic mathematical model, including the key dynamic programming. In Sec. \ref{sec:structural}, structural properties, mainly the  properties of monotonicity under time-varying assumptions, will be discussed. These properties allow us to develop a class of general heuristic policies, we will present the descriptions of the heuristic policy, calculation of its parameters, as well as the outputs of extensive numerical experiments that demonstrate the effectiveness of the heuristic policies in Sec. \ref{sec:heuristics}. Finally, the paper will be concluded in Sec. \ref{sec:conclusions}.

\section{Mathematical Models}
\label{sec:models}

Consider a finite selling horizon of size $N$. At each time period $n=1,2,\ldots, N$, the demand for the product can be characterized by a bivariate random variable $(P_n, Q_n)$, where $Q_n$ denotes the size of the demand, and $P_n$ the price that the buyer is willing to pay. We assume that $(P_n, Q_n)$ takes discrete values in the set,
\begin{align*}
\{p_1, p_2, \ldots,p_I\} \times \{0, 1,2,\ldots, M\}, 
\end{align*}
where $M>0$ is an integer, and $p_i> p_j$ for $i<j$. Denote the discrete probability distribution of the random variable as,
$$\th_{n,i,j}=\pr[P_n=p_i, Q_n=j], \ \ \forall n=1,2,\ldots, N, i=1,2,\ldots, I, j=1,2,\ldots M.$$
We assume that the distributional information on $(P_n, Q_n)$, including the correlations between the demand random variables in different time periods, is completely known to the seller.

At the time period $n=1,2,\ldots, N$, upon seen the realization of the demand, the seller's decision is whether to accept or reject the demand.
If the demand is accepted, then $Q_n$ units will be deducted from the inventory, and the amount of revenue of $P_nQ_n$, will be 
collected. If it is rejected, the demand will not be fulfilled. Note that no partial fulfillment is allowed, therefore, if the remaining inventory is less than $Q_n$, then the demand will be rejected outright. In the 
classical setting of the problem, $(P_n, Q_n)$ are independent. Of course, the seller's goal is to collect maximal expected revenue by time period $N$.

Let us denote $V(n,d)$ as the maximum expected revenue collected starting with $d$ units of inventory at time period $n$. Then, as pointed out in previous literature, see, e.g. \cite{LinLuYao2008}, $V(n,d)$ satisfies the following dynamic programming recursion,
\begin{align}
V(n,d) =&  V(n+1, d)[\Th_{n,0}+\Th_n(d)] \nonumber \\ & \sum_{i=1}^I \sum_{j\le d} \th_{n,i,j}\times 
\max\{p_ij + V(n+1,d-j), V(n+1,d)\}, \label{eqn:dp}
\end{align}
and
\begin{align}
  \label{eqn:dp1}
  V(N,d) =& \ex[\min\{Q_N, d\}P_N],
\end{align}
where
\begin{align}
\label{dp:terminal}
\Th_{n,0}= 1-\sum_{i=1}^I \sum_{j=1}^M  \th_{n,i,j}, \quad \Th_n(d) =\sum_{i=1}^I \sum_{j> d} \th_{n,i,j}.
\end{align}
Once the value function is given, the optimal policy will be naturally checking the difference between the value functions of the different actions, and picking the one that produces the larger return. More specifically, at time $n$ with inventory level at $d$, if the demand realization is $(P_n, Q_n)$, then, we will compare the value of $V(n+1,d)$ which represents the maximum average return from time period $n+1$ until $N$ if we don't accept the request, with the value of $V(n+1, d-Q_n)$,  which represents the maximum average return from time period $n+1$ until $N$ if we do accept the request. If the difference can be offset by the return brought in by the demand, that is $P_nQ_n$, then accept, otherwise reject. Thus, we can see the importance of obtaining the evaluation or approximation of the value function.

\noindent
%{\bf Remark } We will not present a proof for Lemma \ref{lem:monotone},  (i) and (ii) in Lemma \ref{lem:monotone} is, of course, trivial. There exist several different proofs of (iii) and (iv) for different problem settings. (iv) is hinted in \cite{LinLuYao2008}. Meanwhile, the proof that we are going to present covers this basic case.

\section{Structural Properties}
\label{sec:structural}

It has been proved under various settings (discrete or continuous) that when the demand is of unit size, i.e. $Q_n=1$, and the price distributions are independent and identically distributed (IID) then the value function for the stochastic knapsack problem possesses rich first and second order monotonicity properties. More specifically, $V(n,d)$ is increasing and concave in $d$, decreasing and concave in $n$, and is submodular in $(n,d)$. Since the optimal policy is determined quantitatively by $V(n+1,d)-V(n+1,d-1)$, these properties directly lead to monotonic behavior of the optimal policy. 
\begin{itemize}
	\item
	The concavity in $d$ implies that if a certain price class should be accepted with certain amount of inventory on hand, then it should also be accepted if there is less inventory. 
	\item
	The concavity in $n$ means that if one price class should be accepted at certain time, it should also be accepted at any time period after that if the remaining inventory is the same. 
	\item
	Submodularity leads to the fact that the likelihood of one price class being accepted is increasing over time. 
\end{itemize}
In this section, these monotonicity properties are verified under very general assumptions on the price distributions, while we still focus on the unit demand case. Moreover, we will show that the value function, in fact, also satisfies multimodularity, another important monotonicity property that is often instrumental in control problems. Then, in the case that the demand sizes follow general distributions, we show that the value function can be bounded from above and below by value functions for two related problems with unit demand. 

%\subsection{Monotonicity for the Time-Varying Demand}
%\label{sec:monotonicity}

The following theorem establishes the monotonicity properties under a very general assumption of the demand arrivals. 

\begin{thm}
	\label{pro:monotone_time-verying} {\rm Given that the demand size is always one, and $\th_{n, j}$ (short for $\th_{n,1, j}$ ) being a general probability function for time-varying discrete random variable, $V(n,d)$ has the following first and
		second order monotone properties:
		\begin{itemize}
			\item[(i)]
			$V(n,d)$ is decreasing in $n$ and increasing in $d$;
			\item[(ii)]
			$V(n,d)$ is concave in $n$, i.e.,
			%for each $d$,
			$$V(n-1,d)-V(n,d)\le V(n,d)-V(n+1,d);$$
			%holds for $2\le n\le N-1$;
			\item[(iii)]
			$V(n,d)$ is concave in $d$, i.e.,
			%for each $n$,
			$$V(n,d)-V(n,d-1)\ge V(n,d+1)-V(n,d);$$
			%holds for $1\ge d \le W-1$.
			\item[(iv)]
			$V(n,d)$ is submodular in $(n,d)$, i.e.,
			$$V(n,d)-V(n,d-1)\ge V(n+1,d)-V(n+1,d-1);$$
			\item[(v)] A ``submodular-plus'' property
			\begin{eqnarray}
			\label{submodplus} V(n,d+1)-V(n,d)\le V(n+1,d)-V(n+1,d-1).
			\end{eqnarray}
			
		\end{itemize}
	}
\end{thm}

{\startb Proof.}  (i) and (ii) can be easily obtained from the progress of the dynamic programming, similar results and proofs can be found in \cite{talluri2004theory}.

(iii), (iv) and (v) will be established together though induction. 
Note that (v), along with the submodularity in (iv), implies the
concavity in (iii).

First observe that the concavity in $d$, i.e. (iii),  holds
at the boundary, i.e., when $n=N$.

Next, assume that the concavity in $d$ holds for any $n+1$, then, from the definition, we can conclude that,  submodularity holds, i.e.,
$$V(n,d)-V(n,d-1)\ge V(n+1,d)-V(n+1,d-1).$$
Let us establish the inequality in (\ref{submodplus}). Rewrite it as
\begin{eqnarray}
\label{submodplus1} V(n,d+1)-V(n+1,d)\le V(n,d)-V(n+1,d-1).
\end{eqnarray}
The left hand side (LHS) above, upon expanding on $V(n,d+1)$, is
equal to the following (with probability $\th_{n,i,1}$):
\begin{eqnarray*}
	{\rm LHS}&=&p_{i}\ind\{p_{i}+V(n+1,d)\ge V(n+1,d+1)\}\\
	&&+ [V(n+1,d+1)-V(n+1,d)]\ind\{p_{i}+V(n+1,d)< V(n+1,d+1)\}.
\end{eqnarray*}
In the same scenario, the right hand side (RHS) of
(\ref{submodplus1}) becomes:
\begin{eqnarray*}
	{\rm RHS}&=&p_{n}\ind\{p_{i}+V(n+1,d-1)\ge V(n+1,d)\}\\
	&&+ [V(n+1,d)-V(n+1,d-1)]\ind\{p_{i}+V(n+1,d-1)< V(n+1,d)\}.
\end{eqnarray*}
Note that concavity in $d$ (at $n+1$), which has been established
already, leads to
$$p_{i}+V(n+1,d-1) \ge V(n+1,d)\;\imp\;
p_{i}+V(n+1,d) \ge V(n+1,d+1).$$
In this case, ${\rm RHS}={\rm LHS} =p_{i}$.
%also follows from the established concavity.
In the complementary case, i.e., when
$$p_{i}+V(n+1,d-1) < V(n+1,d),$$
we have
\begin{eqnarray*}
	{\rm RHS}= V(n+1,d)-V(n+1,d-1) >p_{i}
\end{eqnarray*}
whereas either ${\rm LHS}=p_{i}$, or
$${\rm LHS}=V(n+1,d+1)-V(n+1,d)\le
V(n+1,d)-V(n+1,d-1)={\rm RHS},$$ due to, again, concavity (in $d$,
at $n+1$).

Having now established the property in (\ref{submodplus}) and the
submodularity in (iv), both for $n$, we can put the two
together to prove the concavity in (iii), for $n$. We can then
repeat this cycle. $\Box$

\noindent
{\bf Remark } Similar results and proofs on concavity and submodularity can also be found in \cite{talluri2004theory}. In general, they are common properties in many systems involving with max and plus operators, see, e.g. Glasserman \& Yao~\cite{glasserman1994monotone}. However,  \eqref{submodplus} does not appear to be common, and hence the way we prove the theorem. Furthermore,  it is closely related to the concept of multimodularity. In addition, property (v) has the following implication of the optimal policy. It says that if a price class is acceptable at certain time, then the same price class is also acceptable next time period, after the inventory has been decreased one unit. 
\begin{defn}
	\label{defn:multiM}
	{\rm A function $f: {\bZ}^n \rightarrow {\mathbb R}$ is called {\it
			multimodular} if the function ${\tilde f}: \bZ^{n+1}\rightarrow
		{\mathbb R}$  defined by
		\begin{align}
		\label{eqn:defn_multiM}
		{\tilde f} (x_0, \bx) = f(x_1-x_0, x_2-x_1, \ldots, x_n-x_{n-1}),
		\quad x_0\in \bZ, \bx\in \bZ^n.
		\end{align}
		is submodular in $n+1$ variable.}
\end{defn}
From this definition, it can be easily verified that (iii), (iv) and (v) imply that the value function $V(t,d)$ is a multimodularity function of $(t,d)$.

As we have seen, the value function for unit demand stochastic knapsack problem is  equipped with rich structural properties, namely, the first and second order monotonicities. It makes them conceptually and computationally more desirable quantities for the purpose of both theoretical study and practical policy implementation. Meanwhile, the structural properties, especially the second order monotonicity, appear to be fairly fragile, and the dynamic programming can not be guaranteed for such properties even when the demand size is drawn from some really common and simple distributions. In the following example, we consider a system with two price classes, and order sizes of $\{1,2,3,4\}$ with a uniform joint distribution. From the three dimensional illustration of the value function (Figure 1.), we can see that convexity can not be expected. 

\bsk
\begin{figure}
	\begin{center}
		\includegraphics[width=4in]{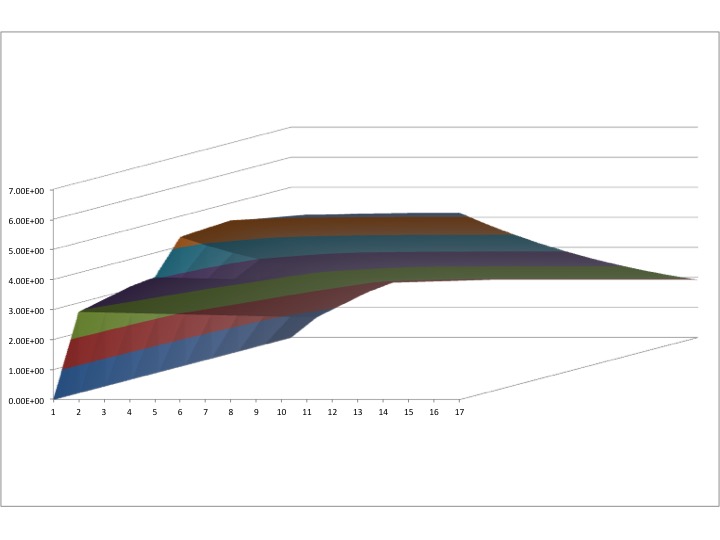}
	\end{center}
	\caption{Value function for a simple example.}
\end{figure}
\bsk

It is therefore natural to seek ways of using the value function of unit demand or other simple functions to approximate or bound the general value function.  From \cite{LinLuYao2008}, we can see that various effective and easy-to-compute lower bounds can be easily obtained, and some of them are even very close to the value function asymptotically, hence we will focus on obtaining an upper bound for the value function.

\begin{thm}
	{\rm For a stochastic knapsack problem as defined in the model section with value function $V(n,d)$, for $n = 1,2,\ldots N$, there exists a unit demand problem with time horizon $NM$ with value function $v(m,d)$, and a mapping $\al: \{1,2, \ldots, N\} \rightarrow \{ 1,2,\ldots, M\}$, such that the following relationship always holds,
		\begin{align*}
		V(n,d) \le v(\al(n), d).
		\end{align*}
	}
\end{thm}
{\bf Proof }
First, it is easy to see that, for any $d\ge 0$, we have, $V(N,d ) \le v((N-1)m+1, d)$, where $v((N-1)M+1, d)$ represents the value function for a system with unit demand, starting from time period $(N-1)M+1$ to time period $NM$, and let the probability for price class $i, i=1,2,\ldots, I$ being $\sum_{m=1}^M\th_{N,i,m}/M$. Therefore, in the original system, the reward is for selling $d$ units in one period, but in the new system, we collect the reward by selling them one by one over $M$ periods. 

In general, any feasible policy at time $n$ with inventory $d$, for the general problem can be expressed as a function $f_{n,d}(p,q)$ taking values in $\{0,1\}$, $f_{n,d}(p,q)=0$ means reject demand with price $p$ and quantity $q$, $f_{n,d}(p,q)=1$ means accept. Thus, we can compute the average revenue garnered in this period, that is $$\sum_{pq} pq\th_{n,p,q} f_{n,d}(p,q).$$
Now, let us consider a unit demand system, from time periods $Mn+1$ to $M(n+1)$. At each time period, the probability for the $i^{th}$ price class will be $p_{i,t}$, and we follow the acceptance/rejection of probability $f_{n,d}(p,q)$, then, the average revenue obtained in this unit demand system will be the same as the one in the general system. The optimal revenue will certainly be an upper bound for $V(n,d)$. $\Box$

\section{Asymptotically Optimal Heuristic Policies}
\label{sec:heuristics}

In \cite{LinLuYao2008}, a family of closed-loop switch-over policies are proposed and studied. Each of the policies identifies a sequence of time epochs. Before the first time epoch, only the demand with the highest price will be accepted. At each subsequent time epoch, one more class of demand, the one with the highest price among those that have not been accepted,  will become acceptable, until at the end all the demand classes will be accepted. The switch-over epochs are purely parameterized by the remaining time and are the results of a nonlinear optimization problem for maximizing average revenue among all switch-over policies of this type. The policies are shown to be asymptotically optimal as the system size grows large, and their efficiency is demonstrated  through extensive numerical experiments. Meanwhile, these policies' emphasis on time certainly limits their strength, especially when the variability of demand size is large and changes over time. In this section, we will extend the idea of closed-loop approach, and propose a family of  two-dimensional switching curve policies, in which thresholds will be determined by both remaining time and inventory, thus enables us to capture the trade-off between time and inventory more effectively. 

Our discussion will first focus on the case of two price classes, then extend to the general cases. In the case of two price classes, the heuristic policy is parameterized by a pair $(t_1, w_1)$. At the beginning, only the high price class is accepted, until either time epoch $t_1$ is reached, or the inventory reaches the level of $w_1+ \sum_{t=t_1+1}^N\ex[Q_{12}(t)]$, then, both classes will be accepted afterwards. Note that $Q_{12}(t)$ represents the random demand for the superposition of the two demand classes  (similarly, $Q_i(t)$ represents the demand for class  $i$), and $ \sum_{t=t_1+1}^N\ex[Q_{12}(t)]$ represents the cumulative average demand starting from time $t_1$ till the end of the time horizon.  In the general case of $M$ classes, the following procedure will be used to produce the switch-over thresholds.
\begin{itemize}
	\item
	Divide classes into two groups, the first group consists of class 1, and the second group consist of the rest of the class.
	\item
	Obtain the optimal $t_1, w_1$ for these two groups. 
	\item
	Consider a new problem for the second group, assuming that time starts at $t_1$, and inventory starts at $W- \sum_{t=1}^{t_1}\ex[Q_1(t)]$, again divide them into two groups, the first group will only have class 2, and the second class will have the rest. Continue with this procedure. 
\end{itemize} 
The output of the procedure will be a sequence of pairs $(t_i, w_i)$, $i=1,2,\ldots, M$, and the heuristic policies: before time $t_i$ or the inventory threshold $w_i + \sum_{t=t_1+1}^N\ex[Q_{i-}(t)] $ can be reached, only classes $1, \ldots, i$ can be accepted, where $Q_{i-}(t)$ represents the superposition demand for classes $1,2,\ldots, i$.  

\noindent
{\bf Remark} 
Our two-dimensional switch-over policies are extensions of the switch-over policies in Lin, Lu and Yao \cite{LinLuYao2008}, hence the asymptotic optimality results obtained in \cite{LinLuYao2008} also apply. It is also worth mentioning that this family of policies are in the similar spirit as those proposed in Bitran and Mondschein \cite{BitranMondschein95}. However, two-dimensional switch-over policies are parameterized, and the selected parameters would likely lead to optimal revenue through optimization over all possible parameters.

%\subsection{Asymptotic Optimality of the Switch-Over Policies}
%\label{sec:asym_opt}

\bsk
%\resizebox{120mm}{!}{\includegraphics*[0in,0in][10in,8in]{revlog_figure3.eps}}
\bsk

%\begin{figure}[!htb]
%\centering
%\includegraphics[scale=.7]{revlog_figure2.eps}
%\caption{Digraph.}
%\label{fig:digraph}
%\end{figure}

%We also like to show that in the case of stationary arrival, and
%linear rewards, the policy that guarantees the maximal reward
%among switch-over policy is asymptotic optimal in certain sense.
%To see that, for any initial condition $T$ and $W$, consider a
%policy that only accept the high priced demand, therefore, its
%average reward will be $p_1\ex[W\wedge
%\sum_{i=1}^{N_1(T)}Q_i^1]^+$. Since $N(T)\rightarrow \infty, a.s.$
%as $T\rightarrow \infty $, therefore, for any $\epsilon>0$, and
%fixed $W$, there exists a $T(W,\epsilon)$, such that for any
%$T>T(W,\epsilon)$, we have
%$$1-\epsilon< \frac{p_1\ex[W\wedge \sum_{i=1}^{N_1(T)}Q_i^1]^+}{p_1W}\le
%1$$ Meanwhile the best of the switch over policy will yield more
%reward than $p_1\ex[W\wedge \sum_{i=1}^{N_1(T)}Q_i^1]^+$, and any
%policy will be collect more revenue that $p_1W$. Hence, the ratio
%between the average reward of the best of the switch-over policies
%and optimal policy will follow the same relationship.

%In the following, we show how to convert the heuristics in
%\cite{bitran} to a $2$-dimensional switch-over policy.

\subsection{Determining the Operational Parameters $(t_1, w_1)$ for the Two Price Class Case}
\label{sec:parameter}

In this section, we will discuss a numerical procedure for determining the switch-over curve for the heuristic policies. Note that, in \cite{LinLuYao2008}, the switch-over epochs are calculated so that the expected total reward is maximized. The derivations heavily rely on the exponential assumption one can not expect in a general setting. For example, one will find it difficult to fit a Poisson model for different arrivals with roughly similar averages, but quite different volatilities. It is also difficult to model time-varying demand, while the time inhomogeneous Poisson is an option for a theoretical model, the computational convenience can not be extended in a straightforward fashion. To cope with these problems, we use a diffusion process to model the demand process, then obtain the switch-over curve parameters with the aid of the calculations of its first passage time at various levels. 

A diffusion process approach provides a reasonable balance between modeling strength and mathematical tractability. As a sophisticated model, it allows us to tract the changes of the first and second moments of the underlying process.  Meanwhile, its mathematical tractability enables the calculations of many useful and important quantities. Rigorous diffusion approximation arguments can be carried out, but it will not be the focus of the current paper, interested readers can consult related literature, see, e.g. \cite{MasseyMandelbaum95}.

In the following, we replace the demand processes of the two classes by two independent Brownian motions $X_1(t)$ and $X_2(t)$, with time-varying drifts $\g_1(t)$ and $\g_2(t)$, and variances $\sigma_1^2(t)$ and $\sigma_2^2(t)$, respectively. Hence, 
\begin{align*} dX_i(t) = \g_i(t) dt +\sigma_i dw_i(t), \quad i=1,2, \end{align*}
where $w_i(t), i=1,2$ are independent standard Brownian motions. From basic calculations for Gaussian processes, we know that, $\ex[X_i(t)]= \g_i(t)$ and $\var[X_i(t)]=\sigma_i^2(t)$. For the best approximation results, let $\g_i(t) = \ex[Q_i(t)]$ and $\sigma_i^2(t)= \var[Q_i(t)]$ for $t=1,2,\ldots, N$, then interpolate the rest values with a continuous function.
%$$\var[\sum_{i=1}^{N_i(t)}D_i]= \ex N_i(t)\var[D_i]+ \ex^2[D_i]\var[N(t)], $$
%$$\var[\sigma_i B_i(t)] =\sigma^2 t. $$
Furthermore, denote the superposition of $X_1(t)$ and $X_2(t)$ as $X_{12}(t)$, and we know that it is a diffusion process with drift $\g_{12}(t)=\g_1(t)+\g_2(t)$ and variance $\sigma_{12}^2(t) = \sigma_1^2(t)+ \sigma_2^2(t)$.

To calculate the average revenue under the diffusion setting, we need to identify the distribution of the first passage time for the approximating one-dimensional diffusion processes. Define the first passage time at level $a>0$ as,
%$$G_i(s,a,b) = \inf\{ t\ge s|X_{i}(t) \ge b\}. $$
\begin{align*}
\tau_a := \inf \{ t\ge 0,X(t)\ge a\},
\end{align*}
for a generic diffusion process $X(t)$ with time-varying drift $\g(t)$ and variance $\sigma(t)$, and $X(0)=0$. Following either a basic change of measure argument, see, e.g. \cite{karatzasbook}, or a partial differential equation argument, see, e.g \cite{2011PhyA}, we know that, the probability density function of the $\tau_a$ bears the following form, 
\begin{eqnarray}
\label{density_passage} f_\tau(t) =
\frac{\sigma^{2}[a+\Lambda(t)]}{\sqrt{4\pi \Sigma(t)^3}}\exp \left[ -\frac{
	(a+\Lambda(t))^2}{4\Sigma(t)}\right],
\end{eqnarray}
with $\Lambda(t) =\int_0^t\g(s)ds$ and $\Sigma(t) = \int_0^t
\sigma^2(s)ds$. 

Suppose that we follow the switch-over policy with parameter $(t_1, w_1)$. According to the policy, before time $t_1\wedge \tau_1$ only the first class will be accepted, where 
$$\tau_1= \inf\left\{t\ge 0\Big| X_1(t) \ge W-w_1-\int_t^N\g_{12}(s)ds\right\}.$$
Equivalently, we have,
$$\tau_1= \inf\left\{t\ge 0\Big| X_1(t)-\int_0^t\g_{12}(s)ds \ge W-w_1-\int_0^N\g_{12}(s)ds\right\}.$$
Denote $Y(t)$ as the diffusion process that satisfies the following stochastic differential equation,
\begin{align*}
dY(t) = [\g_1(t)+ \g_{12}(t) ]dt + \sigma_1(t) dW_1(t).
\end{align*}
Then, $\tau_1$ can be treated as the a first passage time for $Y(t)$, i.e.,
$$\tau_1= \inf\left\{t\ge 0\Big| Y(t) \ge W-w_1-\int_0^N\g_{12}(s)ds\right\}.$$
Therefore, the density function of $\tau_1$ has the following form,
\begin{eqnarray}
\label{density_passage_Y} f_Y(t) =f_\tau(t) =
\frac{\sigma^{2}[{\tilde W}+\Lambda(t)]}{\sqrt{4\pi \Sigma(t)^3}}\exp \left[ -\frac{
	({\tilde W}+\Lambda(t))^2}{4\Sigma(t)}\right]
\end{eqnarray}
with 
\begin{align*}
{\tilde W} & = W-w_1-\int_0^N\g_{12}(s)ds, \\
\Lambda(t) &=\int_0^t[\g_1(s)+\g_{12}(s)]ds, \\
\Sigma(t) &= \int_0^t\sigma^2_{1}(s)ds.
\end{align*} 
Hence, the average rewards collected in this interval will be,
\begin{align*}
V_1(t_1, w_1) & = \ex\left[\int_0^{\tau_1\wedge t_1}p_1X_1(t) dt\right]=\int_0^{\infty} \left(\int_0^{t\wedge t_1} p_1\g_1(s)ds \right)f_Y(dt). 
\end{align*}

%\begin{eqnarray*}
%	V_1= (p_1\g_1) \ex[t_1\wedge \tau_{W-w_1}^{X_{1}}]
%\end{eqnarray*}
In the second step, both classes will be accepted before the inventory is depleted. Let us denote the time of depletion as $\tau_2$. It is easy to see that, $\tau_2$ is the passage time at level $W$ for the process $X_{12}(t)$ starting at $X_1(\tau_1)$. Therefore, the total expected revenue collected in this step can be expressed as, 
\begin{align*}
V_2(t_1, w_1) &= \ex\left[\int_{\tau_1\wedge t_1}^{\tau_2\wedge T}X_{12}(t)\right] = \int_0^{\infty} \int_{t_1}^{\infty} \left(\int_{s\wedge t_1}^{t\wedge T} p_{12}\g_{12}(u)du \right) F_{12}(dt)F_Y(dt). 	
\end{align*}
%\begin{align*}
%	V_2 & = (p_1\g_1+ p_2\g_2) \ex[T\wedge \tau_2-t_1\wedge \tau_1]
%	\end{align*}
%where $\tau_2$ is the stopping time that the entire inventory is depleted. 
Where $F_{12}(dt)$ denotes the density function for the passage time $\tau_2$ for process $X_{12}(t)$, and it is easy to see that a similar expression to \eqref{density_passage_Y} can be used.

In summary, to achieve the best performance, i.e. collecting
maximum total revenue, we need to find the best parameters $(t_1, w_1)$, which will be the output of the following optimization problem,
\begin{eqnarray}
\label{optimization} \max_{t_1, w_1} \quad  V_1(t_1, w_1)+V_2(t_1, w_1).
\end{eqnarray}

\subsection{Numerical Experiments}
\label{sec:numericals}

%\subsection{More on the time-varying demand}
%\subsection{ARMA Demand}
%\subsection{The Performance of the Heuristic Policy}
To examine the efficiency of the two-dimensional switch-over policy, extensive numerical experiments that cover a wide spectrum of parameters and setups are designed and conducted. More specifically, there are three groups of experiments. In the first group, there are only two price classes, with the higher price is set at $1$, and the lower price is randomly generated according to uniform distribution in $(0,1)$, and the size distributions of two classes over time, as well as the probability of no arrival, are also randomly generated. With five instances generated for each case, totally $25$ sets of system parameters are generated. Meanwhile, three different time horizons, $T=10, 20$, and $30$ are selected, and the starting inventory is always set to be $20$. For these totally $125$ systems, we compare the average revenue collected following the heuristic policy, through simulation, against the value function. In eight instances, the value function is within the $97.5$th percentile of the simulation estimates, in $75$ instances, the value function is within the $95$th percentile, $31$ cases for the $80$th percentile, and $11$ cases in the $85$th percentile. In another two groups, we increase the number of price classes to four and eight respectively, in Table 1, the number of instances in different percentiles are listed. Overall, these experiments indicate that the heuristics perform reasonably well. Meanwhile, for each of the system, we also produce the revenue under the one-dimensional switch-over policies, and the two-dimensional policies consistently outperform them by $5-10\%$. 

\vskip 0.5cm

\begin{center}
	{
		\begin{tabular}{r|c|c|c|c|}
			\multicolumn{1}{r}{}
			&\multicolumn{1}{c}{ $97.5$th percentile}
			&\multicolumn{1}{c}{$95$th percentile}
			&\multicolumn{1}{c}{$90$th percentile}
			&\multicolumn{1}{c}{$85$th percentile}\\
			\cline{2-5} Group I  & $8$ & $75$ & $31$ &$11$ \\
			\cline{2-5} Group II & $5$ & $53$ & $52$ & $15$\\
			\cline{2-5} Group III & $4$ & $60$ & $38$ & $23$\\
			\cline{2-5}
		\end{tabular}
		\newline
		\newline 
		Table 1. Performance of the switch-over heuristic policies. 
	}
\end{center}

%\vskip 0.5 cm

\section{Conclusions} 
\label{sec:conclusions}
A version of stochastic knapsack problem with time-varying demand is studied. 
First, we are able to demonstrate the monotonicity of the value function for unit demand stochastic knapsack problem under very general assumptions on the pricing distributions. Then, we extended the switch-over heuristic discussed in 
\cite{LinLuYao2008} to a class of two-dimensional switch-over policies, and demonstrate the efficiency of their performance through extensive numerical experiments. This model can serve as a building block for large and complex revenue management problems for networks, and the results obtained here can be instrumental in the study of such systems, which is our ongoing research.

\vskip 0.5cm

\noindent
%\center{REFERENCES}

%\bibliographystyle{elsarticle-num} 

\bibliographystyle{abbrv}
\bibliography{Lu}


%\begin{thebibliography}{00}


%\bibitem{}

%\end{thebibliography}
\end{document}